\documentclass[a4paper,12pt,reqno]{scrartcl}

\makeatletter
\let\@nohyperfootnotemark\@footnotemark
\let\@nohyperfootnotetext\@footnotetext
\def\nohyperfootnote{\@ifnextchar[\@xfootnote{\stepcounter\@mpfn
     \protected@xdef\@thefnmark{\thempfn}%
     \@nohyperfootnotemark\@nohyperfootnotetext}}
\makeatother

\newcommand\blfootnote[1]{%
  \begingroup
  \renewcommand\thefootnote{}\nohyperfootnote{#1}%
  \addtocounter{footnote}{-1}%
  \endgroup
}

\usepackage[british]{babel}

\usepackage{amsmath,amssymb} % showkeys has amsmath hooks

\usepackage[retainorgcmds]{IEEEtrantools}

\usepackage{xifthen}
\usepackage{enumerate}
\usepackage{bbm} % provides a bb 1
\usepackage{mathrsfs} % provides copperplate \mathscr
\usepackage{amsthm,aliascnt}
\usepackage[final,bookmarksopen=true]{hyperref} % should be last
\usepackage[atend]{bookmark} % should be laster

\usepackage{mathabx}

\newenvironment{eqnarr}{\begin{IEEEeqnarray}{rCl}}{\end{IEEEeqnarray}\ignorespacesafterend}
\newenvironment{eqnarr*}{\begin{IEEEeqnarray*}{rCl}}{\end{IEEEeqnarray*}\ignorespacesafterend}
\newcommand{\eqnarrLHS}[1]{\IEEEeqnarraymulticol{3}{l}{#1} \\ \quad}
\newenvironment{IEEEcases}[1][l]{%[number one][letter ell]
  \left\{%
  \begin{IEEEeqnarraybox}[][c]{#1?l}%
    \IEEEstrut
      }{%
    \IEEEstrut
  \end{IEEEeqnarraybox}%
  \right.%
  \ignorespacesafterend
}

\renewcommand{\eqref}[1]{\hyperref[#1]{(\ref*{#1})}}
\def\pref{\eqref} % w/e
\newcommand{\refpref}[2]{\hyperref[#2]{\ref*{#1}(\ref*{#2})}}
\newcommand{\autorefpref}[2]{\hyperref[#2]{\autoref*{#1}(\ref*{#2})}}

\newcommand{\newsstheorem}[2]{
  \newaliascnt{#1}{dummy}
  \newtheorem{#1}[#1]{#2}
  \aliascntresetthe{#1}
  \expandafter\def\csname #1autorefname\endcsname{#2}
}

\theoremstyle{plain}
  \newsstheorem{thrm}{Theorem}
  \newsstheorem{prop}{Proposition}
  \newsstheorem{cor}{Corollary}
  \newsstheorem{lem}{Lemma}
\theoremstyle{definition}
  \newsstheorem{defn}{Definition}
  \newsstheorem{expl}{Example}
  \newsstheorem{assn}{Assumption}
\theoremstyle{remark}
  \newsstheorem{rem}{Remark}

\def\swap#1#2{\let\arwtempa#1\let#1#2\let#2\arwtempa}
\swap\epsilon\varepsilon

\newcommand{\CE}{\Psi}
\newcommand{\stP}{\mathrm{P}}
\newcommand{\stE}{\mathrm{E}}
\newcommand{\LevP}{\mathbb{P}}
\newcommand{\LevE}{\mathbb{E}}
\newcommand{\stPaz}{\stP^\updownarrow}
\newcommand{\Xaz}{X^\updownarrow}

\newcommand{\Raz}{R^\updownarrow}
\newcommand{\xiaz}{\xi^\updownarrow}

\newcommand*{\define}[1]{\emph{#1}}
\newcommand{\skippar}{\medskip\noindent}

\newcommand{\RR}{\mathbb{R}}
\newcommand{\ZZ}{\mathbb{Z}}
\newcommand{\CC}{\mathbb{C}}

\newcommand{\FF}{\mathscr{F}}
\newcommand{\GG}{\mathscr{G}}

\newcommand{\stproca}[1]{\left(#1\right)_{t \ge 0}}
\newcommand{\stproc}[1]{\stproca{#1_t}}

\newcommand{\GGt}{\stproc{\GG}}

\newcommand{\Ind}{\mathbbm{1}}

\newcommand{\leftsub}[2]{{\protect\vphantom{#2}}_{#1}{#2}}
\newcommand{\lrsub}[3]{\leftsub{#1}{#2}_{#3}}

\newcommand{\downto}{\downarrow}

\newcommand{\wh}{\qquad} % \wh = where; \for conflicts with amscls!
\newcommand{\dd}{\mathrm{d}} % for calculus

\DeclareMathOperator{\sgn}{sgn}

\newcommand*{\newvarcmd}[2]{%
  \newcommand*{#1}[2][]{%
    \begingroup % \sizel and \sizer are local:
      \let\varl\left
      \let\varr\right
      \ifthenelse{\isempty{##1}}{%
        \let\sizel\relax
        \let\sizer\relax
      }{%
        \expandafter\let\expandafter\sizel\csname ##1l\endcsname
        \expandafter\let\expandafter\sizer\csname ##1r\endcsname
      }%
      #2%
    \endgroup
  }
}

\newvarcmd{\abs}{\sizel\lvert #2\sizer\rvert}

\newvarcmd{\norm}{\sizel\lVert #2\sizer\rVert}

\newcommand{\arxivref}[2]{\href{http://arxiv.org/abs/#1}{arXiv:#1} [#2]}

\newcommand{\HG}{\lrsub{2}{F}{1}}

\DeclareMathOperator{\oRe}{Re} % operator (typically roman) style
\DeclareMathOperator{\oIm}{Im} % "
\renewcommand{\Re}{\oRe}
\renewcommand{\Im}{\oIm}
\newcommand{\iu}{\mathrm{i}} % imaginary unit

\newcommand\ltR{\xi}
\newcommand\ltRp{\xi^\prime}
\newcommand\al{\alpha}
\newcommand\MTaux{\widetilde{\MT}}
\newcommand\XRp{R^\prime}
\newcommand\XR{R}
\newcommand\hgbe{\beta}
\newcommand\hgga{\gamma}
\newcommand\hgbeh{\hat\beta}
\newcommand\hggah{\hat\gamma}
\let\xhgbe\hgbe
\let\xhgbeh\hgbeh
\let\xhgga\hgga
\let\xhggah\hggah
\newcommand\tze{T_0}
\newcommand\EF{I}
\newcommand\xhgaux{\zeta}
\newcommand\X{X}

\renewcommand{\ltR}{\xi^{R}}
\renewcommand{\ltRp}{\xi^{R^\prime}}

\DeclareMathOperator{\Res}{Res}
\newcommand{\Resat}[2]{\Res(#1: #2)}

\newcommand{\Ttrans}{\mathcal{T}}
\newcommand{\EssT}{\mathcal{E}}

\newcommand{\DoneyC}{\mathcal{C}}

\newcommand{\MT}{\mathcal{M}} % Mellin transform
\newcommand{\aux}{\zeta} % auxiliary HG process

\newcommand{\betah}{\hat{\beta}}
\newcommand{\gammah}{\hat{\gamma}}
\newcommand{\rhoh}{\hat{\rho}}
\newcommand{\kappah}{\hat{\kappa}}

\newcommand{\ltY}{\xi^{Y}}

\newcommand{\Asym}{\mathcal{A}}
\newcommand{\AHG}{\Asym_{\text{HG}}}
\newcommand{\AXHG}{\Asym_{\text{EHG}}}
\newcommand{\Ast}{\Asym_{\text{st}}}
\newcommand{\Abeta}{\AXHG^{\beta}}
\newcommand{\Abetah}{\AXHG^{\betah}}

\newcommand*{\Gfrac}[2]{\frac{\Gamma(#1)}{\Gamma(#2)}}
\newcommand*{\Indic}[1]{\Ind_{\{#1\}}}

\newcommand{\Ghgsymb}{\HG}
\newcommand*{\Ghg}[4]{\Ghgsymb(#1,#2;#3;#4)}

\IEEEeqnarraydefcolsep{9}{4em}

\usepackage[numbers]{natbib}
\usepackage{doi}

\newcommand{\email}[1]{\href{mailto:#1}{\nolinkurl{#1}}}

\title{The extended hypergeometric class of L\'evy processes}
\author{
  A.\ E.\ Kyprianou%
  \footnote{%
    University of Bath, UK. \email{a.kyprianou@bath.ac.uk}
  }
  \and
  J.\ C.\ Pardo%
  \footnote{%
    CIMAT, Mexico. \email{jcpardo@cimat.mx}
  }
  \and
  A.\ R.\ Watson%
  \footnote{%
    University of Z\"urich, Switzerland. \email{alexander.watson@math.uzh.ch}.
    Parts of this work were completed while the third author was at the
    University of Bath, UK, and at CIMAT, Mexico.
  }
}
\date{\today}

\hypersetup{
  pdftitle={The extended hypergeometric class of L\'evy processes},
  pdfauthor={A. E. Kyprianou, J. C. Pardo and A. R. Watson}
}

\begin{document}

\maketitle

\begin{abstract}
  \noindent
  With a view to computing fluctuation identities related to
  stable processes, we review and extend the class of
  hypergeometric L\'evy processes explored in
  \citet{KP-HG}. We give the Wiener--Hopf factorisation
  of a process in the extended class, and characterise its exponential
  functional. Finally, we give three concrete examples
  arising from transformations of stable processes.%
  \blfootnote{%
    \textit{Keywords and phrases:} L\'evy processes, hypergeometric L\'evy processes,
    extended hypergeometric L\'evy processes,
    Wiener--Hopf factorisation, exponential functionals,
    stable processes, path-censored stable process,
    conditioned stable processes, hitting distributions, hitting probabilities.%
  }%
  \blfootnote{%
    \textit{MSC 2000 classification:} 60G51, 60G18, 60G52.%
  }
\end{abstract}

\section{Introduction}

The simple definition of a L\'evy process---a stochastic process with
stationary independent increments---has been sufficient to fuel a vast
field of study for many decades, and L\'evy processes have been employed
in many successful applied models. However, historically there have been
few classes of processes for which many functionals could be computed
explicitly. In recent years, the field
has seen a proliferation of examples which have proved to be more
analytically tractable; in particular, we single out
spectrally negative L\'evy processes
\cite{KKR-scale}, Lamperti-stable processes \cite{CPP-LS,CKP-explicit},
$\beta$- and $\theta$-processes \cite{Kuz-beta,Kuz-theta},
and finally the inspiration for this work, hypergeometric
L\'evy processes \cite{KP-HG,KKPvS,KPR-n-tuple}. It is also worth mentioning
that the close relationship which appears to hold between hypergeometric
L\'evy processes and stable processes has also allowed the computation
of several identities for the latter; see \cite{KP-HG,KPR-n-tuple,KPW-cens}.

In this work, we review the hypergeometric class of L\'evy processes
introduced by \citet{KP-HG}, and introduce a new class
of extended hypergeometric processes which have many similar properties.
In particular, for an extended hypergeometric process $\xi$
we compute the Wiener--Hopf factors
and find that its ladder height processes
are related to Lamperti-stable subordinators; and we are able to
characterise explicitly the distribution of the exponential
functional of $\xi/\delta$ for any $\delta > 0$.
We also give three examples of processes connected via the
Lamperti representation to $\alpha$-stable processes, which fall into
the hypergeometric class when $\alpha \le 1$, and into the extended
hypergeometric class when $\alpha > 1$, and give some
new identities for the stable process when $\al > 1$.

\skippar
We will first discuss the results of \citet{KP-HG}.
For a choice of parameters $(\beta,\gamma,\betah,\gammah)$ from the set
\[ \AHG = \bigl\{ \beta \le 1, \, \gamma \in (0,1), \,
  \betah \ge 0, \, \gammah \in (0,1) \bigr\} , \]
we define
\[ \psi(z)
  = - \Gfrac{1-\beta+\gamma-z}{1-\beta-z}
  \Gfrac{\betah+\gammah+z}{\betah+z} , \]
%for $z \in \CC \setminus \{ 1-\beta+\gamma+k , -(\hat\beta+\hat\gamma+k) : k =0,1,2,\dotsc\}$.
which we view as a meromorphic function on $\CC$.
We say that a L\'evy process $\xi$ is a member of the \define{hypergeometric class} of L\'evy processes
if it has Laplace exponent $\psi$, in the sense that 
\begin{equation}\label{e:ce}
  \LevE[e^{z \xi_1}] = e^{\psi(z)}, \wh z \in \iu\RR.
\end{equation}
Note that, in general, when the Laplace exponent $\psi$ of a
L\'evy process $\xi$ is a meromorphic function, the relation
\eqref{e:ce} actually holds on any neighbourhood of $0 \in \CC$
which does not contain a pole of $\psi$; thus, in this article we will generally
not specify the domain of Laplace exponents which may arise.

In \cite{KP-HG}, it is shown that for any choice of parameters in $\AHG$, there
is a L\'evy process with Laplace exponent $\psi$, and they find its Wiener--Hopf factorisation,
in the following sense.

The \define{(spatial) Wiener--Hopf factorisation} of a L\'evy process $\xi$ with Laplace
exponent $\psi$ consists of the
equation
\[ \psi(z) = - \kappa(-z) \kappah(z), \wh z \in \iu \RR, \]
where $\kappa$ and $\hat\kappa$ are the Laplace exponents of subordinators
$H$ and $\hat H$, respectively, this time in the
sense that $\LevE\bigl[e^{-\lambda H_1}\bigr] = e^{-\kappa(\lambda)}$
for $\Re \lambda \ge 0$.
The subordinators $H$ and $\hat H$ are known as the ascending and descending
ladder heights, and are related via a time-change to the running maximum and running minimum of
the process $\xi$. For more details, we refer the
reader to \cite[Chapter 6]{Kyp}. The insight into the structure of $\xi$
given by the Wiener--Hopf factorisation allows one to simplify
first passage problems for $\xi$; see \cite[Chapter 7]{Kyp} for a
collection of results.

The paper \cite{KP-HG} computes that
\[ \kappa(z) = \frac{\Gamma(1-\beta+\gamma+z)}{\Gamma(1-\beta+z)},
  \qquad
  \hat\kappa(z) = \frac{\Gamma(\hat\beta+\hat\gamma+z)}{\Gamma(\hat\beta+z)},
\]
thus demonstrating that the ascending and descending ladder height processes
% sensitive formatting here
are\linebreak\mbox{Lamperti-stable} subordinators (see \cite{CPP-LS}).

Kuznetsov and Pardo also consider the exponential functional of a hypergeometric
L\'evy process $\xi$.
For each $\delta > 0$, the random variable
% $\hat\beta > 0$, and one defines the exponential functional
\[ I(\xi/\delta) = \int_0^\infty e^{-\xi_t/\delta} \, \dd t \]
is a.s.\ finite provided that $\xi$ drifts to $+\infty$.
This random variable is known as the \define{exponential functional}
of the L\'evy process $\xi$,
and it has been studied extensively in general;
the paper of \citet{BY-ef} gives a survey of the literature,
and mentions, among other aspects, applications to diffusions in
random environments, mathematical finance and fragmentation theory.
In the context of self-similar Markov processes, the exponential functional
appears in the entrance law of a pssMp started at zero (see,
for example, \citet{BY-ent}), and
\citet{Par-upper} relates the exponential functional of a L\'evy process
to envelopes of its associated pssMp; furthermore, it is related
to the hitting time of points for pssMps, and we shall make use of it
in this capacity in our example of \autoref{s:radial}.

For the purpose of characterising the distribution of $I(\xi/\delta)$,
its Mellin transform
\[ \MT(s) = \LevE[I(\xi/\delta)^{s-1}] , \]
is useful. For $\xi$ in the hypergeometric class, $\MT$ was calculated by
\cite{KP-HG} in terms of gamma and double gamma functions;
we will recall and extend this \autoref{s:ef}.

\skippar
We now give a brief outline of the main body of the paper.
In section 2, we demonstrate that the parameter set $\AHG$ may be extended
by changing the domains of the two parameters $\beta$ and $\betah$, and find the Wiener--Hopf
factorisation of a process $\xi$ in this new class, identifying explicitly
the ladder height processes. In section 3, we find an expression for the
Mellin transform $\MT$ in this new case, making use of an auxiliary
hypergeometric L\'evy process. In section 4, we cover three examples where
the extended hypergeometric class is of use, on the way
extending the result of \cite{CPP-radial} on the Wiener--Hopf
factorisation of the Lamperti representation associated with
the radial part of a stable process.

\section{The extended hypergeometric class}
\label{s:HG}
\label{s:XHG}

We begin by defining the set of admissible parameters
\[ \AXHG = \bigl\{ \beta \in [1,2], \,
  \gamma,\hat\gamma \in (0,1), \,
  \hat\beta \in [-1,0]; \,
  1-\beta+\betah+\gamma \ge 0, \,
  1-\beta+\betah+\gammah \ge 0 \bigr\} . \]
We are interested in proving the existence of, and investigating the properties of,
a L\'evy
process $\xi$ whose Laplace exponent is given by the meromorphic function
\[ 
  \psi(z)
  = - \Gfrac{1-\beta+\gamma-z}{1-\beta-z}
  \Gfrac{\betah+\gammah+z}{\betah+z},
  \wh z \in \CC,
\]
when $(\beta,\gamma,\betah,\gammah) \in \AXHG$.
  
To allow for more concise expressions below, we also define
\[ \eta = 1-\beta+\gamma+\betah+\gammah . \]
We now give our main result on the existence and properties of $\xi$.

\begin{prop}\label{t:WHF XHG}
  There exists a L\'evy process $\xi$ such that $\LevE[e^{z \xi_1}] = e^{\psi(z)}$.
  Its Wiener-Hopf factorisation may be expressed as
  \[ \psi(z)
    = - (-\betah-z) \Gfrac{1-\beta+\gamma-z}{2-\beta-z}
    \times (\beta-1+z)
    \Gfrac{\betah+\gammah+z}{1+\betah+z} . \]
  Its L\'evy measure possesses the density
  \begin{equation}\label{e:XHGLM}
    \pi(x) = \begin{cases}
      - \dfrac{\Gamma(\eta)}
      {\Gamma(\eta-\gammah) \Gamma(-\gamma)}
      e^{-(1-\beta+\gamma)x}
      \Ghg{1+\gamma}{\eta}{\eta-\gammah}{e^{-x}},
      & x > 0, \\
      - \dfrac{\Gamma(\eta)}{\Gamma(\eta-\gamma)\Gamma(-\gammah)}
      e^{(\betah+\gammah)x}
      \Ghg{1+\gammah}{\eta}{\eta-\gamma}{e^{x}},
      & x < 0 ,
    \end{cases}
  \end{equation}
  where $\Ghgsymb$ is the Gauss hypergeometric function.

  If $\beta \in (1,2)$ and $\betah \in (-1,0)$, the process $\xi$ is killed at rate
  \[ q = \Gfrac{1-\beta+\gamma}{1-\beta} \Gfrac{\betah+\gammah}{\betah} . \]
  Otherwise, the process has infinite lifetime and:
  \begin{enumerate}[(i)]
  \item\label{i:WH:+} $\xi$  drifts to $+\infty$ if
    $\beta > 1, \, \betah = 0$.
  \item\label{i:WH:-} $\xi$ drifts to $-\infty$ if
    $\beta = 1, \, \betah < 0$.
  \item\label{i:WH:osc} $\xi$ oscillates if
    $\beta = 1, \, \betah = 0$.
    In this case, $\xi$ is a hypergeometric L\'evy process.
  \end{enumerate}

  Furthermore, the process $\xi$ has no Gaussian component,
  and is of bounded variation with zero drift when $\gamma+\gammah < 1$
  and of unbounded variation when $\gamma+\gammah \ge 1$.

\begin{proof}
We remark that there is nothing to do in case
\pref{i:WH:osc} since such processes are analysed in \cite{KP-HG};
however, the proof we give below also carries through in this case.

We will first identify the proposed ascending and descending ladder
processes. Once we have shown that $\psi$ really is the
Laplace exponent of a L\'evy process, this will be the proof
of the Wiener-Hopf factorisation.

\newcommand{\lesbd}{\upsilon}
\newcommand{\leasc}{\upsilon}
\newcommand{\ledesc}{\hat\upsilon}
Before we begin, we must review the definitions of special subordinators
and the $\Ttrans$-transformations of subordinators. Suppose that $\lesbd$ is
the Laplace exponent of a subordinator $H$, in the sense that
$\LevE[e^{-zH_1}] = e^{-\lesbd(z)}$. $H$ is said to be a
\define{special subordinator}, and $\lesbd$ a \define{special
Bernstein function}, if the function
\[ \lesbd^*(z) = z/\lesbd(z), \wh z \ge 0, \]
is also the Laplace exponent of a subordinator. The function $\lesbd^*$ is
said to be \define{conjugate} to $\lesbd$.
Special Bernstein functions
play an important role in potential theory; see, for example,
\cite{SV-special} for more details.

Again taking $\lesbd$ to be the Laplace exponent of a subordinator,
not necessarily special, we define, for $c \ge 0$, the transformation
\[
  \Ttrans_{c}\lesbd(z)
  = \frac{z}{z+c}
  \lesbd(z+c),
  \wh z \ge 0.
\]
It is then known (see \cite{Gne-rrcs,KP-CT}) that
$\Ttrans_c\lesbd$ is the Laplace exponent of a subordinator.
Furthermore, if $\lesbd$ is in fact a special Bernstein function,
then $\Ttrans_c\lesbd$ is also a special Bernstein function.
% ;
% this is not difficult to show, but a proof may be found in
% \cite[Proposition 5.2]{KPW-cens}.

We are now in a position to identify the ladder height processes
in the Wiener--Hopf factorisation of $\xi$.
Let the (proposed) ascending factor be given, for $z \ge 0$, by
\[ \kappa(z) = (-\betah+z) \Gfrac{1-\beta+\gamma+z}{2-\beta+z} . \]
Then some simple algebraic manipulation shows that
\[ \kappa(z) = \bigl(\Ttrans_{-\betah} \leasc\bigr)^*(z) , \]
with
\[ \leasc(z) = \frac{\Gamma(2-\beta+\betah+z)}
    {\Gamma(1-\beta+\betah+\gamma+z)} , \]
provided that $\leasc$ is a special Bernstein function.
This follows immediately from Example 2 of \citet{KR-scale},
% To verify this,
% we apply Example 2 of \citet{KR-scale}. This immediately gives
% the result
under the constraint $1-\beta+\betah+\gamma \ge 0$ which
is included in the parameter set $\AXHG$.
Furthermore, we note that $\leasc$ is in fact the Laplace exponent
of a Lamperti-stable subordinator (see \cite{CPP-LS}), although we will
not use this fact.

%% Removed this but it is true:
% Furthermore, we may even identify this $\leasc$ as the Laplace
% exponent of a Lamperti-stable subordinator with parameters
% \[ \bigl(q, \, \LSa, \, \LSb, \, \LSc, \, \LSd \bigr)
%   = \biggl( \Gfrac{2-\beta+\betah}{1-\beta+\betah+\gamma}, \,
%     1-\gamma, \,
%     \beta-\betah-\gamma, \,
%     -\frac{1}{\Gamma(\gamma-1)}, \,
%     0 \biggr) ,
% \]
% notation as in \cite{KPW-cens}.

Proceeding similarly for the descending factor, we obtain
\[ \kappah(z) = (\beta-1+z)\Gfrac{\betah+\gammah+z}{1+\betah+z}
  = \bigl( \Ttrans_{\beta-1} \ledesc \bigr)^*(z),
  \wh z \ge 0,\]
where
\[ \ledesc(z) = \Gfrac{2-\beta+\betah+z}{1-\beta+\betah+\gammah+z} , \wh z \ge 0, \]
and again the function $\ledesc$ is a special Bernstein function provided
that $1-\beta+\betah+\gammah \ge 0$.
As before, $\ledesc$ is the Laplace exponent of a Lamperti-stable
subordinator.
%% As above:
% Again, $\ledesc$ is the Laplace
% exponent of a Lamperti-stable subordinator with parameters
% \[ \bigl(q, \, \LSa, \, \beta, \, c, \, \LSd \bigr)
%   = \biggl( \Gfrac{2-\beta+\betah}{1-\beta+\betah+\gammah}, \,
%     1-\gammah, \,
%     \beta-\betah-\gammah, \,
%     -\frac{1}{\Gamma(\gammah-1)}, \,
%     0 \biggr) .\]

We have now shown that both $\kappa$ and $\kappah$ are Laplace
exponents of subordinators; we wish to show that the function
\begin{equation*}\label{e:friends}
  \psi(z) = -\kappa(-z) \kappah(z)
\end{equation*}
is the Laplace exponent of a L\'evy process. For this purpose
we apply the theory of \define{philanthropy} developed by
\citet[chapter 7]{Vig-thesis}. This states, in part, that
it is sufficient for both of the subordinators corresponding to
$\kappa$ and $\kappah$ to be `philanthropists', which means
that their L\'evy measures possess decreasing densities.

We recall our discussion of $\Ttrans$-transforms and special
Bernstein functions. We have already stated that when $\lesbd$ is
a special Bernstein function, then so is
$\Ttrans_c\lesbd$; furthermore, one may show that
its conjugate satisfies
\[
  (\Ttrans_c\lesbd)^*(z)
  = \EssT_c\lesbd^*(z) + \lesbd^*(c), \wh z \ge 0,
\]
where $\EssT_c$ is the \define{Esscher transform}, given by
\[ \EssT_c\lesbd^*(z) = \lesbd^*(z+c)-\lesbd^*(c) , \wh z \ge 0. \]
The Esscher transform of the Laplace exponent of any subordinator is again
the Laplace exponent of a subordinator; and
if the subordinator corresponding to $\lesbd^*$ possesses
a L\'evy density $\pi_{\lesbd^*}$, then the L\'evy density of
$\EssT_c \lesbd^*$ is given by $x \mapsto e^{-cx} \pi_{\lesbd^*}(x)$,
for $x > 0$.

Returning to our Wiener-Hopf factors, we have
\[\kappa(z) = \bigl(\Ttrans_{-\betah} \leasc\bigr)^*(z)
  = \EssT_{-\betah}\leasc^*(z) + \leasc^*(-\betah),
  \wh z \ge 0, \]
where $\leasc^*$ is the
Laplace exponent conjugate to $\leasc$.
Now, $\leasc$ is precisely the type of special Bernstein function
considered in \cite[Example 2]{KR-scale}. In that work, the
authors even establish that the subordinator corresponding to
$\leasc^*$ has a decreasing L\'evy density $\pi_{\leasc^*}$.
Finally, the L\'evy density of the subordinator corresponding
to $\kappa$ is $x \mapsto e^{\betah x} \pi_{\leasc^*}(x)$,
and this is then clearly also decreasing.

Hence, we have shown that the subordinator whose Laplace exponent
is $\kappa$ is a philanthropist. By a very similar argument,
the subordinator corresponding to $\kappah$ is also a philanthropist.
As we have stated, the theory developed by Vigon now shows that
the function $\psi$ really is the Laplace exponent of a L\'evy process
$\xi$, with the Wiener--Hopf factorisation claimed.

We now proceed to calculate the L\'evy measure of $\xi$.
A fairly simple way to do this is to make use of the theory of
`meromorphic L\'evy processes', as developed in
\citet{KKP-mero}.
We first show that $\xi$ is in the meromorphic class.
Initially suppose that
\begin{equation}\label{e:gt assum}%\tag{\textit{gt}}
  1-\beta+\betah+\gamma > 0, \;
  1-\beta+\betah+\gammah > 0;
\end{equation}
we will relax this assumption later.
Looking at the expression for $\psi$, we see that it has zeroes
$(\zeta_n)_{n \ge 1}, (-\hat\zeta_n)_{n \ge 1}$
and (simple) poles $(\rho_n)_{n \ge 1}, (-\hat\rho_n)_{n \ge 1}$
given as follows:
\begin{IEEEeqnarray*}{rCl9rCl"l}
  \zeta_1 &=& -\betah,
  & \zeta_n &=& n -\beta , &n \ge 2, \\
  &&&\rho_n &=& n-\beta+\gamma, &n \ge 1, \\
  \hat\zeta_1 &=& \beta-1,
  &\hat\zeta_n &=& \betah+n-1, &n \ge 2, \\
  &&&\hat\rho_n &=& \betah+\gammah+n-1, &n \ge 1,
\end{IEEEeqnarray*}
and that they satisfy the interlacing condition
\[
  \dotsb < - \hat\rho_2 < - \hat\zeta_2 < -\hat\rho_1 < -\hat\zeta_1
  < 0 < \zeta_1 < \rho_1 < \zeta_2 < \rho_2 < \dotsb .
\]
To show that $\xi$ belongs to the meromorphic class, one applies
\cite[Theorem 1(v)]{KKP-mero} when $\xi$ is killed, and
\cite[Corollary 2]{KKP-mero} in the unkilled case. The proof is
a routine calculation using the Weierstrass representation \cite[8.322]{GR}
to expand $\kappa$ and $\kappah$ as infinite products, and we
omit it for the sake of brevity.

We now calculate the L\'evy density.
For a process in the meromorphic class, it is known that
the L\'evy measure has a density of the form
\begin{equation}\label{e:LD of XHG}
  \pi(x) = \Indic{x > 0} \sum_{n \ge 1} a_n \rho_n e^{-\rho_n x}
  + \Indic{x < 0} \sum_{n \ge 1} \hat a_n \hat\rho_n e^{\hat\rho_n x},
\end{equation}
for some coefficients $(a_n)_{n \ge 1}$, $(\hat a_n)_{n \ge 1}$,
where the $\rho_n$ and $\hat\rho_n$ are as above.
Furthermore, from \cite[equation (8)]{KKP-mero}, we see that
\[ a_n \rho_n = -\Resat{\psi(z)}{z = \rho_n} , \]
and correspondingly for $\hat a_n \hat \rho_n$. (This remark
is made on p.~1111 of \cite{KKP-mero}.) From here it is simple
to compute
\[ a_n \rho_n
  = - \frac{(-1)^{n-1}}{(n-1)!}
  \frac{1}{\Gamma(1-\gamma-n)}
  \frac{\Gamma(\eta+n-1)}{\Gamma(\eta-\gammah+ n - 1)},
  \wh n \ge 1,
\]
and similarly for $\hat a_n \hat \rho_n$.
The expression \eqref{e:XHGLM}
follows by substituting in \eqref{e:LD of XHG} and using
the series definition of the hypergeometric function.

Thus far we have been working under the assumption that
\eqref{e:gt assum} holds. Suppose now that this fails and we have,
say, $1-\beta+\betah+\gammah=0$. Then $\zeta_1=\rho_1$,
which is to say the first zero-pole pair to the right of
the origin is removed.
It is clear that $\xi$ still falls into the meromorphic
class, and indeed, our expression for $\pi$ remains valid:
although the initial pole $\rho_1$ no longer exists, the
corresponding coefficient $a_1\rho_1$ vanishes as well.
Similarly, we may allow $1-\beta+\betah+\gamma=0$,
in which case the zero-pole pair to the
left of the origin is removed; or we may allow both
expressions to be zero, in which case both pairs are removed.
The proof carries through in all cases.

The claim about the large time behaviour of $\xi$ follows from the
Wiener-Hopf factorisation: $\kappa(0) = 0$ if and only if the
range of $\xi$ is a.s.\ unbounded above, and $\kappah(0) = 0$
if and only if the range of $\xi$ is a.s.\ unbounded below;
so we need only examine the values of $\kappa(0),\,\kappah(0)$
in each of the four parameter regimes.

Finally, we prove the claims about the Gaussian
component and variation of $\xi$. This
proof proceeds along the same lines as that in \cite{KP-HG}.
Firstly, we observe
using \cite[formula 8.328.1]{GR} that
\begin{equation}\label{e:psi asym}
  \psi(\iu\theta)
  = O\bigl(\abs{\theta}^{\gamma+\gammah}\bigr), \wh \abs{\theta} \to \infty.
\end{equation}
Applying \cite[Proposition I.2(i)]{Ber-Levy} shows that $\xi$ has no Gaussian component.
Then, using \cite[formulas 9.131.1 and 9.122.2]{GR}, one
sees that
\[ \pi(x) = O\bigl(\abs{x}^{-(1+\gamma+\gammah)}\bigr), \wh x \to 0, \]
and together with the necessary and sufficient condition
$\int_\RR (1\wedge\abs{x}) \pi(x)\,\dd x < \infty$
for bounded variation, this proves the claim
about the variation. In the bounded variation case,
applying \cite[Proposition I.2(ii)]{Ber-Levy} with \eqref{e:psi asym}
shows that $\xi$ has zero drift.
\end{proof}
\end{prop}

We propose to call this the \define{extended hypergeometric class}
of L\'evy processes.

\begin{rem}
If $\xi$ is a process in the extended hypergeometric class,
with parameters
$(\beta,\gamma,\betah,\gammah)$, then the dual process $-\xi$ also
lies in this class, and has parameters\linebreak
$(1-\betah,\gammah,1-\beta,\gamma)$.
\end{rem}

\begin{rem}
  We remark here that one may instead extend the parameter range $\AHG$
  by moving only $\beta$, or only $\betah$. To be precise, both
  \[ \Abeta =
    \bigl\{ \beta \in [1,2], \, \gamma,\gammah \in (0,1), \,
    \betah \ge 0; \, 1-\beta+\betah+\gamma \le 0, \,
    1-\beta+\betah+\gammah \ge 0 \bigr\} \]
  and
  \[ \Abetah = 
    \bigl\{ \beta \le 1, \, \gamma,\gammah \in (0,1), \,
    \betah \in [-1,0]; \, 1-\beta+\betah+\gamma \ge 0, \,
    1-\beta+\betah+\gammah \le 0 \bigr\} \]
  are suitable parameter regimes, and
  one may develop a similar theory for such processes; for instance,
  for parameters in $\Abeta$, one has the Wiener--Hopf factors
  \[ \kappa(z) = \frac{\Gamma(2-\beta+\gamma+z)}{\Gamma(2-\beta+z)},
    \qquad
    \kappah(z) = \frac{\beta-1+z}{\beta-1-\gamma+z}
    \frac{\Gamma(\betah+\gammah+z)}{\Gamma(\betah+z)} .
  \]
  However, we are not aware of any examples of processes in these classes.
\end{rem}

\section{The exponential functional}
\label{s:ef}

\newcommand{\psieta}{\tilde{\psi}}
\newcommand{\thetaaux}{\tilde{\theta}}

Suppose that $\xi$ is a L\'evy process in the extended
hypergeometric class with $\beta > 1$, which is to say 
either $\xi$ is killed or it drifts to $+\infty$.

We are then
interested in the \define{exponential functional} of the
process, given by
\[ I(\xi/\delta) = \int_0^\infty e^{-\xi_t/\delta} \, \dd t , \]
for $\delta >0$. (Since $\xi/\delta$ is not in the extended hypergeometric
class, we are studying exponential functionals of a slightly larger
collection of processes.) This is an a.s.\ finite random variable under
the conditions we have just outlined.

It will emerge that the best way to characterise the distribution
of $I(\xi/\delta)$ is via its \define{Mellin transform},
\[ \MT(s) = \LevE\bigl[I(\xi/\delta)^{s-1}\bigr] , \]
whose domain of definition will be a vertical strip in the complex
plane to be determined.

In the
case of a hypergeometric L\'evy process with $\betah > 0$, it was shown in \cite{KP-HG}
that the Mellin transform of the exponential functional is given by
\[
  \MT_{\text{HG}}(s) = C \Gamma(s) 
  \frac{G((1-\beta)\delta+s;\delta)}{G((1-\beta+\gamma)\delta + s;\delta)}
  \frac{G((\hat\beta+\hat\gamma)\delta + 1-s;\delta)}{G(\hat\beta\delta+1-s;\delta)} ,
\]
holds
for $\Re s \in (0,1+\hat\beta\delta)$, where $C$ is a normalising constant
such that $\MT_{\text{HG}}(1)=1$,
and $G$ is the double gamma function; see \cite{KP-HG} for a definition
of this special function.

Our goal in this section is
the following result, which characterises the law of the exponential functional for the extended
hypergeometric class.

\begin{prop}\label{t:MT XHG}
Suppose that $\xi$ is a L\'evy process in the extended hypergeometric
class with $\beta > 1$.
Define $\theta = \delta(\beta-1)$.

Then, the Mellin transform $\MT$ of $I(\xi/\delta)$
is given by
\begin{equation} \MT(s) = 
  c \MTaux(s)
  \Gfrac{\delta(1-\beta+\gamma)+s}{-\delta\betah+s}
  \Gfrac{\delta(\beta-1)+1-s}{\delta(\betah+\gammah)+1-s} , 
  \wh \Re s \in (0, 1+\theta),
  \label{e:MT of XHG}
\end{equation}
where $\MTaux$ is the Mellin transform of $I(\aux/\delta)$,
and $\aux$ is an auxiliary L\'evy process in the hypergeometric
class, with parameters $(\beta-1,\gamma,\betah+1,\gammah)$. The constant
$c$ is such that $\MT(1) = 1$.
\begin{proof}
The process $\xi/\delta$ has Laplace exponent $\psi_\delta$ given by
$\psi_\delta(z) = \psi(z/\delta)$. The relationship with $\aux$ arises
from the following calculation:
\begin{eqnarr}
  \psi_\delta(z)
  &=& \frac{-\betah-z/\delta}{1-\beta+\gamma-z/\delta}
    \frac{\beta-1+z/\delta}{\betah+\gammah+z/\delta}
    \Gfrac{2-\beta+\gamma-z/\delta}{2-\beta-z/\delta}
    \Gfrac{1+\betah+\gammah+z/\delta}{1+\betah+z/\delta} \nonumber \\
  &=& \frac{-\betah-z/\delta}{1-\beta+\gamma-z/\delta}
    \frac{\beta-1+z/\delta}{\betah+\gammah+z/\delta}
    \psieta_\delta(z), \label{e:xi and eta}
\end{eqnarr}
where $\psieta_\delta$ is the Laplace exponent of a L\'evy process $\aux/\delta$,
with $\aux$ as in the statement of the theorem.

Denote by $f(s)$ the right-hand side of \eqref{e:MT of XHG}.
The proof now proceeds via the `verification result'
\cite[Proposition 2]{KP-HG}.

Recall that a L\'evy process with Laplace exponent $\phi$ is said to
satisfy the \define{Cram\'er condition with Cram\'er number $\theta$} if
there exists $z_0 < 0$ and $\theta \in (0,-z_0)$
such that $\phi(z)$ is defined for all $z \in (z_0,0)$ and
$\phi(-\theta) = 0$.

Inspecting the Laplace exponent $\psi_\delta$ reveals that $\xi/\delta$
satisfies the Cram\'er condition with Cram\'er number
$\theta = \delta(\beta - 1)$.

Furthermore, $\aux/\delta$ satisfies the Cram\'er condition
with Cram\'er number $\thetaaux = \delta(\betah+1)$.
It follows from \cite[Lemma 2]{Riv-re2} that $\MTaux(s)$
is finite in the strip
$\Re s \in (0, 1+\thetaaux)$; and by the properties of
Mellin transforms of positive random variables,
it is analytic and zero-free in its domain of definition.
The constraints in the parameter set $\AXHG$
ensure that $\thetaaux \ge \theta$; this, together with
inspecting the right-hand
side of \eqref{e:MT of XHG} and comparing again with the
conditions in $\AXHG$, demonstrates that $\MT(s)$
is analytic and zero-free in the strip
$\Re s \in (0,1+\theta)$.

We must then check the
functional equation $f(s+1) = -s f(s)/ \psi_\delta(-s)$,
for $s \in (0,\theta)$.
Apply \eqref{e:xi and eta} to write
\begin{eqnarr*}
  -\frac{s}{\psi_\delta(-s)}
  &=& - \frac{s}{\psieta_\delta(-s)}
    \frac{1-\beta+\gamma+s/\delta}{-\betah+s/\delta}
    \frac{\betah+\gammah-s/\delta}{\beta-1-s/\delta} \\
  &=& \frac{\MTaux(s+1)}{\MTaux(s)}
    \frac{\delta(1-\beta+\gamma)+s}{-\delta\betah+s}
    \frac{\delta(\betah+\gammah)-s}{\delta(\beta-1)-s} \\
  &=& \frac{\MTaux(s+1)}{\MTaux( s)}
    \Gfrac{-\delta\betah+s}{-\delta\betah+s+1}
    \Gfrac{\delta(1-\beta+\gamma)+s+1}{\delta(1-\beta+\gamma)+s} \\
  && {} \times \Gfrac{\delta(\betah+\gammah)+1-s}{\delta(\betah+\gammah)-s}
    \Gfrac{\delta(\beta-1)-s}{\delta(\beta-1)+1-s},
\end{eqnarr*}
making use of the same functional equation for the Mellin transform
$\MTaux$. It is then clear that the right-hand side is equal to $f(s+1)/f(s)$.

Finally, it remains to check that
$\abs{f(s)}^{-1} = o(\exp(2\pi\abs{\Im(s)}))$, as
$\abs{\Im s} \to \infty$, uniformly in $\Re s \in (0,1+\theta)$.
%% Abbreviated this to the below:
%By the formula
%\cite[8.328.1]{GR},
%\[ \abs{\Gamma(x+\iu y)}
%  \sim \sqrt{2\pi}
%  e^{-\frac{\pi}{2}\abs{y}}
%  \abs{y}^{\frac{1}{2}-x},
%  \for \abs{y} \to \infty, \]
%we find that
%\[
%  \log \abs{\Gfrac{\delta(1-\beta+\gamma)+s}{-\delta\betah+s}
%  \Gfrac{\delta(\beta-1)+1-s}{\delta(\betah+\gammah)+1-s}}^{-1}
%  \sim
%  \delta(\gamma-\gammah)\log\abs{\rIm s} = o(\abs{\rIm s}),
%\]
The following asymptotic relation may be derived from Stirling's asymptotic
formula for the gamma function:
\begin{equation}\label{e:asym}
  \log \Gamma(z) = z \log z - z + O(\log z),
\end{equation}
and since Stirling's asymptotic formula is uniform in
$\abs{\arg(z)} < \pi-\omega$ for any choice of $\omega > 0$,
it follows that \eqref{e:asym} holds uniformly
in the strip $\Re s \in (0,1+\theta)$; see \cite[Chapter 8, \S 4]{Olv-asym}.
We thus obtain
\[ 
  \log \abs[bigg]{\Gfrac{\delta(1-\beta+\gamma)+s}{-\delta\betah+s}
  \Gfrac{\delta(\beta-1)+1-s}{\delta(\betah+\gammah)+1-s}}^{-1}
  = O( \log s) = o(\Im s), \]
and comparing this with the proof of \cite[Theorem 2]{KP-HG}, where
the asymptotic behaviour of $\MTaux(s)$ is given, we see that
this is sufficient for our purposes.

Hence, $\MT(s) = f(s)$ when $\Re s \in (0,1+\theta)$.
\end{proof}
\end{prop}

This Mellin transform may be inverted to give an expression for the density
of $I(\xi/\delta)$ in terms of series whose terms are defined iteratively,
but we do not pursue this here. For details of this approach, see
\cite[\S 4]{KP-HG}.

\section{Three examples}
\label{s:ex}

It is well-known that hypergeometric L\'evy processes appear as
the Lamperti transforms of stable processes killed passing below
zero, conditioned to stay positive and conditioned to
hit zero continuously; see \cite[Theorem 1]{KP-HG}.
In this section we briefly present three additional examples in which
the extended hypergeometric class comes into play.
The examples may all be obtained in the same way: we begin with a stable process,
modify its path in some way to obtain a positive, self-similar Markov process,
and then apply the Lamperti transform to obtain a new L\'evy process.
We therefore start with a short description of these concepts.

% These paras are from the potentials paper and the path-censored paper.
\skippar
\phantomsection\label{s:stable}%
We work with the (strictly) stable process $X$ with scaling parameter $\alpha$ and
positivity parameter $\rho$, which is defined as follows. For $(\alpha,\rho)$
in the set
\begin{eqnarr*}
  \Ast
  &=&
  \{ (\alpha,\rho) : \alpha \in (0,1), \, \rho \in (0,1) \}
  \cup \{ (\alpha,\rho) = (1,1/2) \}\\
  &&
  {} \cup \{ (\alpha,\rho) : \alpha \in (1,2), \, \rho \in (1-1/\alpha, 1/\alpha) \},
\end{eqnarr*}
and let $X$, with probability laws $(\stP_x)_{x \in \RR}$,
be the L\'evy process with characteristic exponent
\[ \CE(\theta)
  = \begin{cases}
    c\abs{\theta}^\alpha (1  - \iu\beta\tan\frac{\pi\alpha}{2}\sgn\theta)
    & \alpha \in (0,2) \setminus \{1\}, \\
    c\abs{\theta} & \alpha = 1,
  \end{cases}
  \wh \theta \in \RR,
\]
where $c = \cos(\pi\alpha(\rho-1/2))$ and
$\beta = \tan(\pi\alpha(\rho-1/2))/\tan(\pi\alpha/2)$; by this we mean that
$\stE_0[e^{\iu\theta X_1}] = e^{-\CE(\theta)}$.
This L\'evy process has absolutely continuous L\'evy measure with density
\newcommand{\cpex}{\frac{\Gamma(\alpha+1)}{\Gamma(\alpha\rho)\Gamma(1-\alpha\rho)}}
\newcommand{\cmex}{\frac{\Gamma(\alpha+1)}{\Gamma(\alpha\rhoh)\Gamma(1-\alpha\rhoh)}}
\[ c_+ x^{-(\alpha+1)} \Indic{x > 0}
  + c_- \abs{x}^{-(\alpha+1)}\Indic{x < 0},
  \wh x \in \RR , \]
where
\[ c_+ = \cpex, \qquad c_- = \cmex \]
and $\rhoh = 1-\rho$.

The parameter set $\Ast$ and the characteristic exponent $\CE$ represent,
up a multiplicative constant in $\Psi$, all (strictly) stable processes
which jump in both directions, except for
Brownian motion and the symmetric Cauchy processes with non-zero drift.

The choice of $\alpha$ and $\rho$ as parameters is explained as follows.
$X$ satisfies the \define{$\alpha$-scaling property}, that
\begin{equation}\label{e:scaling}
  \text{under } \stP_x \text{, the law of }
  (cX_{t c^{-\alpha}})_{t \ge 0} \text{ is } \stP_{cx} ,
\end{equation}
for all $x \in \RR, \, c > 0$. The second parameter
satisfies $\rho = \stP_0(X_t > 0)$.

\skippar
A \define{positive self-similar Markov process} (\define{pssMp}) with
\define{self-similarity index} $\alpha > 0$ is a standard Markov process
$Y = (Y_t)_{t\geq 0}$ with filtration $\GGt$ and probability laws
$(\stP_x)_{x > 0}$, on $[0,\infty)$, which has $0$ as an absorbing state and
which satisfies the scaling property \eqref{e:scaling} (with $Y$ in place of $X$).
Here, we mean ``standard'' in the sense of \cite{BG-mppt},
which is to say, $\GGt$ is a complete, right-continuous filtration,
and $Y$ has c\`adl\`ag paths and is strong Markov
and quasi-left-continuous.

\skippar
\phantomsection
\label{s:Lamperti}%
In the seminal paper \cite{Lam-ssLT}, Lamperti describes a one to one correspondence
between pssMps and L\'evy processes, which we now outline.
It may be worth noting that we have presented a slightly
different definition of pssMp from Lamperti; for the connection, see
\cite[\S 0]{VA-Ito}.

Let
$S(t) = \int_0^t (Y_u)^{-\alpha}\, \dd u .$
This process is continuous and strictly increasing until $Y$ reaches zero.
Let $(T(s))_{s \ge 0}$ be its inverse, and define
\[ \xi_s = \log Y_{T(s)} \qquad s\geq 0. 
\]
Then $\xi : = (\xi_s)_{s\geq 0}$ is a L\'evy process started at $\log x$, possibly killed at an independent
exponential time; the law of the L\'evy process and the rate of killing do not depend
on the value of $x$. The real-valued process $\xi$ with probability laws
$(\LevP_y)_{y \in \RR}$ is called the
\define{L\'evy process associated to $Y$}, or the \define{Lamperti transform of $Y$}.

An equivalent definition of $S$ and $T$, in terms of $\xi$ instead
of $Y$, is given by taking
$T(s) = \int_0^s \exp(\alpha \xi_u)\, \dd u$
and $S$ as its inverse. Then,
\begin{equation}
  \label{Lamp repr}
  Y_t = \exp(\xi_{S(t)}) 
\end{equation}
for all $t\geq 0$, and this shows that the Lamperti transform is a bijection.

%% True but cut for brevity:
%Let $T_0 = \inf\{ t > 0: Y_t = 0 \}$ be the first hitting time of the absorbing state zero. Then the large-time behaviour of $\xi$ can be described by the behaviour of $Y$ at $T_0$, as follows:
%\begin{enumerate}[(i)]
%  \item If $T_0 = \infty$ a.s., then $\xi$ is unkilled and either oscillates or drifts to $+ \infty$.
%  \item If $T_0 < \infty$ and $Y_{T_0 -} = 0$ a.s., then $\xi$ is unkilled and drifts to $-\infty$.
%  \item If $T_0 < \infty$ and $Y_{T_0 -} > 0$ a.s., then $\xi$ is killed.
%\end{enumerate}
%It is proved in \cite{Lam-ssLT} that the events mentioned above satisfy a zero-one law independently of $x$, and so the three possibilites above are an exhaustive classification of pssMps.

\subsection{The path-censored stable process}
\label{s:cens}

\newcommand{\Ain}{\gamma}

Let $X$ be the stable process defined in \autoref{s:stable}.
In \cite{KPW-cens}, the present authors considered a `path-censored'
version of the stable process, formed by erasing the time spent
in the negative half-line. To be precise, define
\[ A_t = \int_0^t \Indic{X_s > 0} \, \dd s, \wh t \ge 0, \]
and let $\Ain(t) = \inf\{ s \ge 0 : A_s > t \}$ be its right-continuous
inverse. Also define
\[ T_0 = \inf\{t \ge 0 : X_{\gamma(t)} = 0 \}, \]
which is finite or infinite a.s.\ accordingly as $\alpha > 1$ or
$\alpha \le 1$. Then, the process
\[ Y_t = X_{\Ain(t)} \Indic{t < T_0}, \wh t \ge 0, \]
is a pssMp, called the \define{path-censored stable process}.

In Theorems 5.3 and 5.5 of \cite{KPW-cens}, it is shown that
the
Laplace exponent $\psi^Y$
% characteristic exponent $\Psi$
of the Lamperti transform
$\ltY$ associated to $Y$ is given by
\begin{equation*}
  \psi^Y(z)
  = \frac{\Gamma(\alpha\rho-z)}{\Gamma(-z)}
  \frac{\Gamma(1-\alpha\rho+z)}{\Gamma(1-\alpha+z)},
%   \Psi(\theta)
%   = \frac{\Gamma(\alpha\rho-\iu\theta)}{\Gamma(-\iu\theta)}
%   \frac{\Gamma(1-\alpha\rho+\iu\theta)}
%     {\Gamma(1-\alpha+\iu\theta)} , \wh \theta \in \RR,
\end{equation*}
and there it is remarked that when $\alpha \le 1$, this
process is in the hypergeometric class with parameters
\[ \bigl(\beta,\gamma,\hat\beta,\hat\gamma\bigr)
  = \bigl(1, \alpha\rho, 1-\alpha, \alpha\rhoh \bigr). \]
It is readily seen from our definition
that, when $\alpha > 1$, the process $\ltY$
is in the extended hypergeometric class, with the same
set of parameters.

From the Lamperti transform we know that
\[ I(-\alpha\ltY)
  = \inf\{ u \ge 0 : Y_u = 0 \}
  = \int_0^{T_0} \Indic{X_t > 0}\, \dd t,
\]
where the latter is the occupation time of $(0,\infty)$ up to first
hitting zero for the stable process. This motivates the following
proposition, whose proof is a direct application of \autoref{t:MT XHG}.

\begin{prop}
  The Mellin transform of the random variable $I(-\alpha\ltY)$
  is given, for $\Re s \in (\rho-1/\alpha,2-1/\alpha)$, by
  \[
    \MT(s)
    = c \frac{G(2/\alpha-1+s; 1/\alpha)}{G(2/\alpha-\rho+s; 1/\alpha)}
    \frac{G(1/\alpha+\rho+1-s; 1/\alpha)}{G(1/\alpha+1-s; 1/\alpha)}
    \frac{\Gamma(1/\alpha-\rho+s)}{\Gamma(\rho+1-s)}
    \Gamma(2-1/\alpha-s),
  \]
  where $c$ is a normalising constant such that $\MT(1)=1$.
\end{prop}

\begin{rem}
  When $X$ is in the class $\DoneyC_{k,l}$ introduced by \citet{Don-Ckl}, which is to say
  \[ \rho + k = l/\alpha,  \]
  for $k,l \in \ZZ$,
  equivalent expressions in terms of gamma and trigonometric functions
  may be found
  via repeated application of certain identities
  of the double-gamma function; see, for example,
  \cite[equations (19) and (20)]{KP-HG}.
  %\citet[equation (4.6)]{Kuz-extr}.

  For example, when $k,l \ge 0$, one has
  \begin{eqnarr*}
    \MT(s)
    &=& c
    (-1)^{l}
    (2\pi)^{l(1/\alpha-1)}
    (1/\alpha)^{l(1-2/\alpha)}
    \frac{\Gamma(\frac{1-l}{\alpha}+k+s)}{\Gamma(\frac{l}{\alpha}+1-k-s)}
    \frac{\Gamma(2-1/\alpha-s)}{\Gamma(2-l-\alpha-\alpha s)} \\
    && {} \times
    \prod_{j=1}^l
    \Gamma(j/\alpha+1-s)\Gamma(2/\alpha-(j/\alpha+1-s))
    \prod_{i=0}^{k-1}
    \frac{\sin(\pi\alpha(s+i))}{\pi} ,
  \end{eqnarr*}
  and when $k < 0, \, l \ge 0$,
  \begin{eqnarr*}
    \MT(s)
    &=& c
    (-1)^{l}
    (2\pi)^{l(1/\alpha-1)}
    (1/\alpha)^{l(1-2/\alpha)} \\
    && {} \times
    \frac{\Gamma(\frac{1-l}{\alpha}+k+s)\Gamma(2-1/\alpha-s)
        \Gamma(l+1+\alpha-\alpha s)\Gamma(2-l+\alpha k + \alpha s)}
      {\Gamma(\frac{l}{\alpha}+1-k-s)} \\
    && {} \times 
    \prod_{j=1}^l
    \Gamma(j/\alpha+1-s)\Gamma(2/\alpha-(j/\alpha+1-s))
    \prod_{i=2}^{-k-1}
    \frac{\pi}{\sin(\pi\alpha(s-i))} .
  \end{eqnarr*}
  Similar expressions may be obtained when $k \ge 0, \, l < 0$ and $k,l < 0$.
\end{rem}

\subsection{The radial part of the symmetric stable process}
\label{s:radial}

If $X$ is a symmetric stable process---that is, $\rho=1/2$---then
the process
\[ R_t = \abs{X_t}, \wh t \ge 0, \]
is a pssMp, which we call the \define{radial part} of $X$.
The Lamperti transform, $\ltR$, of this process was studied by
\citet{CPP-radial} in dimension $d$; these authors
computed the Wiener--Hopf factorisation of $\ltR$ under the assumption
$\alpha < d$, finding that the process is a hypergeometric L\'evy process.
Using the extended hypergeometric class,
we extend this result, in one dimension, by finding the
Wiener--Hopf factorisation when $\alpha > 1$.

In \citet{KKPW-T0}, the following theorem is proved using the work of
\citet{CPP-radial}.

%% For this section (up to the proof of the Mellin transform)
%% there are a lot of custom macros because it is imported
%% from my thesis. They are defined in defs-thesis.tex.
%% Sorry about this.
\begin{thrm}[Laplace exponent]\label{t:ltR CE}
  The
  Laplace exponent
  % characteristic exponent
  of the L\'evy process $2\ltR$
  is given by
  \begin{equation}\label{e:ltR CE}
    \psi^R(2 z) = -2^{\alpha}
    \frac{\Gamma(\alpha/2-z)}{\Gamma(-z)}
    \frac{\Gamma(1/2+z)}{\Gamma((1-\alpha)/2+z)}.
%     \CER(2\theta) = 2^{\al}
%     \frac{\Gamma(\al/2-\iu\theta)}{\Gamma(-\iu\theta)}
%     \frac{\Gamma(1/2+\iu\theta)}{\Gamma((1-\al)/2+\iu\theta)},
%     \wh \theta \in \RR.
  \end{equation}
\end{thrm}

We now identify the Wiener--Hopf factorisation of $\ltR$, which will
depend on the value of $\al$. However, note the factor
$2^{\al}$ in \eqref{e:ltR CE}. In the context of the
Wiener--Hopf factorisation, we could ignore this
factor by picking an appropriate normalisation of local time; however,
another approach is as follows.

\newcommand\LERp{\psi^\prime}
Let us write $\XRp = \frac{1}{2}\XR$,
and denote by $\ltRp$ the Lamperti transform of $\XRp$.
Then the scaling of space on the level of the self-similar process
is converted by the Lamperti transform into a scaling of time, so
that $\ltR_s = \log 2 + \ltRp_{s 2^{\al}}$. In particular,
if we write $\LERp$ for the characteristic exponent of
$\ltRp$, it follows that $\LERp = 2^{-\al}\psi^R$.
This allows us to disregard the inconvenient constant factor
in \eqref{e:ltR CE}, if we work with $\ltRp$ instead
of $\ltR$.

The following corollary is now simple when we bear in mind
the hypergeometric class of L\'evy processes introduced in
\autoref{s:HG}. We emphasise that this Wiener--Hopf factorisation
was derived by different methods in \cite[Theorem 7]{CPP-radial}
for $\al < 1$, though not $\al = 1$.

\def\templab{Wiener--Hopf factorisation,
  \texorpdfstring{$\alpha \in (0,1]$}{alpha in (0,1]}}
\begin{cor}[\templab]\label{c:R WHF low al}
  The Wiener--Hopf factorisation of $2\ltRp$ when $\al \in (0,1]$
  is given by
  \[
    \LERp(2z)
    = - \frac{\Gamma(\alpha/2-z)}{\Gamma(-z)}
    \times
    \frac{\Gamma(1/2+z)}{\Gamma(1-\alpha)/2+z)}
  \]

% \[
%   \CERp(2\theta) =
%   \frac{\Gamma(\al/2-\iu\theta)}{\Gamma(-\iu\theta)}
%   \times
%   \frac{\Gamma(1/2+\iu\theta)}{\Gamma((1-\al)/2+\iu\theta)},
%   \wh \theta \in \RR,
% \]
  and 2$\ltRp$ is a L\'evy process of the hypergeometric class
  with parameters
  \[ (\hgbe,\hgga,\hgbeh,\hggah) = (1, \al/2, (1-\al)/2, \al/2). \]
%The ladder
%processes $H, \hat H$ associated with $2\xi$ are as follows.
%$H$ is a Lamperti-stable subordinator (see \cite{CPP-LS},
%and \cite[\S 5.2]{KPW-cens} for notation)
%with
%\[ (q,\LSa,\beta,\LSd,c)
%  = \biggl(0,\alpha/2,1,0,-\frac{2^\alpha}{\Gamma(-\alpha/2)}\biggr) \]
%and $\hat H$ is a Lamperti-stable subordinator with
%\[ (q,\LSa,\beta,\LSd,c)
%  = \biggl(\frac{\Gamma(1/2)}{\Gamma((1-\alpha)/2)},
%  \alpha/2,(1+\alpha)/2,0,-\frac{1}{\Gamma(-\alpha/2)}\biggr). \]
\begin{proof}
It suffices to compare the characteristic exponent with that
of a hypergeometric L\'evy process.
\end{proof}
\end{cor}

When $\al > 1$, the process $\ltRp$ is not a hypergeometric L\'evy process;
however, it is in the extended hypergeometric class, and we therefore
have the following result, which is new.

\begin{thrm}[Wiener--Hopf factorisation,
  \texorpdfstring{$\alpha \in (1,2)$}{alpha in (1,2)}]%
\label{t:R WHF high al}
The Wiener--Hopf factorisation of $2\ltRp$ when $\al \in (1,2)$
is given by
\begin{equation}\label{e:ltR WHF high index}
  \LERp(2z)
  =
  -\biggl(\frac{\alpha-1}{2}-z\biggr)
  \frac{\Gamma(\alpha/2-z)}{\Gamma(1-z)}
  \times
  z
  \frac{\Gamma(1/2+z)}{\Gamma((3-\alpha)/2+z)}
% \CERp(2\theta)
%   = 
%     \biggl(\frac{\al-1}{2} - \iu\theta\biggr)
%     \frac{\Gamma(\al/2 - \iu\theta)}{\Gamma(1-\iu\theta)}
%     \times
%     \iu\theta
%     \frac{\Gamma(1/2+\iu\theta)}{\Gamma((3-\al)/2+\iu\theta)},
\end{equation}
and $2\ltRp$ is a L\'evy process in the extended hypergeometric
class, with parameters
\[ (\xhgbe,\xhgga,\xhgbeh,\xhggah)
  = (1, \al/2, (1-\al)/2, \al/2) . \]
  
\begin{proof}
Simply use \autoref{t:ltR CE};
using the formula $x\Gamma(x) = \Gamma(x+1)$
yields \eqref{e:ltR WHF high index}. That this is indeed the Wiener-Hopf
factorisation follows once we recognise $2\ltRp$ as a process
in the extended hypergeometric class, and apply
\autoref{t:WHF XHG}.
\end{proof}
\end{thrm}

\skippar
As an illustration of the utility of the extended hypergeometric class,
we will now derive an expression for the Mellin transform of the 
exponential functional for the dual process $-\ltRp$. This quantity
is linked by the Lamperti representation to the hitting time of zero
for $X$; see \autoref{s:Lamperti}. In particular, if
\[ \tze = \inf\{t \ge 0 : \X_t = 0 \}, \]
we have that
\begin{equation}\label{e:T0 EF}
  \tze
  = \int_0^\infty e^{\al\ltR_t} \, \dd t
  = \int_0^\infty e^{\al\ltRp_{2^{\al}t}} \, \dd t
  = 2^{-\al} \int_0^\infty e^{\al \ltRp_{s}} \, \dd s
  = 2^{-\al} \EF(-\al \ltRp)
  .
\end{equation}

Since $-2\ltRp$ is an extended hypergeometric L\'evy process
with parameters
$\bigl(\frac{\al+1}{2}, \frac{\al}{2}, 0, \frac{\al}{2} \bigr)$,
which drifts to $+\infty$, we may apply the theory just developed to
compute the Mellin transform of $I(-\al \ltRp)$.
Denote this by $\MT$; that is,
\[ \MT(s)
  = \LevE\bigl[\EF(-\al\ltRp)^{s-1}\bigr], \]
for some range of $s \in \CC$ to be determined.

\begin{prop}\label{p:MT symm}
For $\Re s \in (-1/\al, 2-1/\al)$,
\begin{eqnarr}
  \stE_1[T_0^{s-1}]
  &=&
  2^{-\alpha(s-1)}
  \MT(s) \nonumber \\
  &=& 
  2^{-\alpha(s-1)}
  \frac{\sqrt{\pi}}{\Gamma(\frac{1}{\al})\Gamma(1-\frac{1}{\al})}
  \frac{\Gamma(1+\frac{\al}{2}-\frac{\al s}{2})}
    {\Gamma(\frac{1-\al}{2} + \frac{\al s}{2})}
  \Gamma(\tfrac{1}{\al}-1+s)
  \frac{\Gamma(2-\frac{1}{\al}-s)}{\Gamma(2-s)}.
  \label{e:radial:MT}
\end{eqnarr}
\begin{proof}

Let $\xhgaux$ be a hypergeometric L\'evy process with parameters
$\bigl(\frac{\al-1}{2}, \frac{\al}{2}, 1, \frac{\al}{2} \bigr)$,
and denote by $\MTaux$ the Mellin transform of the exponential
functional $I(\al/2\cdot\xhgaux)$, which is known to be finite for
$\Re s \in (0,1+2/\al)$ by the argument in the proof of \autoref{t:MT XHG}.

We can then use \autoref{t:MT XHG} to make the following calculation,
provided that \mbox{$\Re s \in (0,2-1/\al)$}.
Here $G$ is the double gamma function, as defined in \cite[S\S 3]{KP-HG},
and we use \cite[equation (25)]{KP-HG} in the third line
and the identity $x\Gamma(x) = \Gamma(x+1)$ in the final line.
For normalisation constants $C$ (and $C^\prime$) to be determined, we have
%\begin{equation}\label{e:MT inter1}
\[
\begin{IEEEeqnarraybox}[][c]{rCl}
  \MT(s)
  &=& C \MTaux(s)
    \frac{\Gamma(\frac{1}{\al} + s)}{\Gamma(s)}
    \frac{\Gamma(2-\frac{1}{\al}-s)}{\Gamma(2-s)} \\
  &=& C \frac{G\bigl( \frac{3}{\al} - 1 + s ; \frac{2}{\al}\bigr)}
      {G\bigl(\frac{3}{\al}+s; \frac{2}{\al} \bigr)}
    \frac{G\bigl(\frac{2}{\al}+2-s; \frac{2}{\al} \bigr)}
      {G\bigl(\frac{2}{\al}+1-s; \frac{2}{\al} \bigr)}
    \Gamma(s)
    \frac{\Gamma(\frac{1}{\al} + s)}{\Gamma(s)}
    \frac{\Gamma(2-\frac{1}{\al}-s)}{\Gamma(2-s)} \\
  &=& C \frac{\Gamma\bigl(1+\frac{\al}{2}-\frac{\al}{2}s\bigr)}
      {\Gamma\bigl(\frac{3-\al}{2} + \frac{\al}{2}s \bigr)}
    \Gamma(\tfrac{1}{\al} + s)
    \frac{\Gamma(2-\frac{1}{\alpha}-s)}{\Gamma(2-s)} \\
  &=& C^\prime
    \frac{\Gamma(1+\frac{\al}{2}-\frac{\al s}{2})}
    {\Gamma(\frac{1-\al}{2} + \frac{\al s}{2})}
    \Gamma(\tfrac{1}{\al}-1+s)
    \frac{\Gamma(2-\frac{1}{\al}-s)}{\Gamma(2-s)}.
  %\yesnumber\label{e:radial:MT}
\end{IEEEeqnarraybox}
\]
%\end{equation}
The condition $\MT(1) = 1$ means that we can calculate
\begin{equation}\label{e:MT Cp}
  C^\prime
  = \frac{\sqrt{\pi}}{\Gamma(\frac{1}{\al})\Gamma(1-\frac{1}{\al})} ,
\end{equation}
and this gives the Mellin transform explicitly,
for $\Re s \in (0,2-1/\al)$.

We now expand the domain of $\MT$.
Note that, in contrast to the general case of \autoref{t:MT XHG}, the
right-hand side of \eqref{e:radial:MT} is well-defined when
\mbox{$\Re s \in (-1/\al, 2-1/\al)$}, and is indeed analytic in
this region.
(The reason for this difference is
the cancellation of a simple pole and zero at the point $0$.)
Theorem 2 of \cite{LS-acf} shows that, if the Mellin transform
of a probability measure is analytic in a neighbourhood of the
point $1 \in \CC$, then it is analytic in a strip
$\Re s \in (a,b)$, where $-\infty \le a < 1 < b \le \infty$;
and futhermore, the function has singularities at $a$ and $b$,
if they are finite.
It then follows that the right-hand side of \eqref{e:radial:MT}
must actually be equal to $\MT$ in all of $\Re s \in (-1/\al, 2-1/\al)$,
and this completes the proof.
\end{proof}
\end{prop}

We remark that the distribution of $T_0$ has been characterised
previously by \citet{YYY-laws} and \citet{Cor-thesis}, using
rather different methods; and
the Mellin transform above was also obtained,
again via the Lamperti transform but without the extended hypergeometric class,
in \citet{KKPW-T0}.

\skippar
It is also fairly straightforward to produce the following
hitting distribution. Define
\[ \sigma_{-1}^1 = \inf\{ t \ge 0 : X_t \notin [-1,1]\}, \]
the \define{first exit time of $[-1,1]$} for $X$.
We give the
distribution of the position of the symmetric stable
process $X$ at time $\sigma_{-1}^1$,
provided this occurs before $X$ hits zero.
Note that when $\alpha \in (0,1]$, the
process does not hit zero, so the distribution is simply that found
by \citet{Rog-hit}.

\begin{prop}
  Let $X$ be the symmetric stable process with $\alpha \in (1,2)$.
  Then, for $\abs{x} < 1$, $y > 1$,
  \begin{eqnarr*}
    \eqnarrLHS{
      \stP_x\bigl(\abs[big]{X_{\sigma_{-1}^1}} \in \dd y; \; \sigma_{-1}^1 < T_0\bigr)/\dd y
    }
    &=&
    \frac{\sin(\pi\alpha/2)}{\pi}
    \abs{x}
    (1-\abs{x})^{\alpha/2}
    y^{-1}
    (y-1)^{-\alpha/2}
    (y-\abs{x})^{-1}
    \\
    && {} + 
    \frac{1}{2}
    \frac{\sin(\pi\alpha/2)}{\pi}
    y^{-1}
    (y-1)^{-\alpha/2}
    \abs{x}^{(\alpha-1)/2}
    \int_0^{1-\abs{x}}
    t^{\alpha/2-1}
    (1-t)^{-(\alpha-1)/2}
    \, \dd t.
  \end{eqnarr*}
  \begin{proof}
    The starting point of the proof is the `second factorisation identity'
    \cite[Exercise 6.7]{Kyp},
    \[
      \int_0^\infty e^{-qz} \LevE\Bigl[e^{-\beta(\xi_{S_z^+}-z)} ; \; S_z^+ < \infty\Bigr] \, \dd x
      = \frac{\kappa(q) - \kappa(\beta)}{(q-\beta)\kappa(q)},
      \wh q,\beta > 0,
    \]
    where
    \[ S_z^+ = \inf\{ t \ge 0 : \xi_t > z \}. \]
    We now invert in $q$ and $z$, in that order; this is a lengthy but routine calculation,
    and we omit it.
    We then apply the Lamperti transform: if
    $g(z,\cdot)$ is the density of the measure $\LevP(\xi_{S_z^+} - z \in \cdot;\;S_z^+ < \infty)$, then
    \[ \stP_x\bigl(\abs[big]{X_{\sigma_{-1}^1}} \in \dd y; \; \sigma_{-1}^1 < T_0\bigr)
      = y^{-1} g(\log\abs{x}^{-1},\log y) ,
    \]
    and this completes the proof.
  \end{proof}
\end{prop}

The following hitting probability emerges after integrating
in the above proposition.

\begin{cor}
  For $\abs{x}<1$,
  \[
    \stP_x(T_0<\sigma_{-1}^1)
    =
    (1-\abs{x})^{\alpha/2}
    - \frac{1}{2} \abs{x}^{(\alpha-1)/2}
    \int_0^{1-\abs{x}}
    t^{\alpha/2-1}
    (1-t)^{-(\alpha-1)/2}
    \, \dd t .
  \]
\end{cor}

Finally, it is not difficult to produce the following slightly more general result.
Applying the Markov property at time $T_0$ gives
\begin{eqnarr*}
  \stP_x\bigl( X_{\sigma_{-1}^1} \in \dd y; \; \sigma_{-1}^1\bigr)
  &=& \stP_x(X_{\sigma_{-1}^1} \in \dd y)
  - \stP_x(X_{\sigma_{-1}^1} \in \dd y; \; T_0 < \sigma_{-1}^1) \\
  &=& \stP_x(X_{\sigma_{-1}^1} \in \dd y)
  - \stP_x(T_0 < \sigma_{-1}^1) \stP_0(X_{\sigma_{-1}^1} \in \dd y).
\end{eqnarr*}
The hitting distributions on the right-hand side were found by
\citet{Rog-hit}, and substituting yields the following corollary.

\begin{cor}
  For $\abs{x}<1$, $\abs{y}>1$,
  \begin{IEEEeqnarray*}{rCll}
    \IEEEeqnarraymulticol{4}{l}{
      \stP_x(X_{\sigma_{-1}^1} \in \dd y; \; \sigma_{-1}^1<T_0)/\dd y
    } \\
    \quad &=&
    \IEEEeqnarraymulticol{2}{l}{
      \frac{\sin(\pi\alpha/2)}{\pi}
      (1-x)^{\alpha/2}
      (1+x)^{\alpha/2}
      (y-1)^{-\alpha/2}
      (y+1)^{-\alpha/2}
      (y-x)^{-1}
    } \\
    && {} - {} &
    \biggl[
    (1-\abs{x})^{\alpha/2}
    - \frac{1}{2} \abs{x}^{(\alpha-1)/2}
    \int_0^{1-\abs{x}}
    t^{\alpha/2-1}
    (1-t)^{-(\alpha-1)/2}
    \, \dd t
    \biggr]
    \\
    &&& {} \times
    \frac{\sin(\pi\alpha/2)}{\pi}
    (y-1)^{-\alpha/2}
    (y+1)^{-\alpha/2}
    y^{-1} .
  \end{IEEEeqnarray*}
\end{cor}

\subsection{The radial part of the symmetric stable process conditioned to avoid zero}
\label{s:radial az}

Above, we computed the Lamperti transform
of the pssMp $R^\prime = \frac{1}{2}\abs{X}$, where $X$ was a symmetric
stable process, and called it $\ltRp$. In this section we consider
instead the symmetric stable process conditioned to avoid zero,
and obtain its Lamperti transform.

In \cite{Pan-cond}, \citeauthor{Pan-cond} shows (among many other results)
that the function
\[ h(x) =
  \begin{IEEEcases} % IEEEcases adds some more space
  {-} \Gamma(1-\alpha)
  \frac{\sin (\pi\alpha\rhoh)}{\pi}
  x^{\alpha-1} , & x > 0, \\
  {-} \Gamma(1-\alpha)
  \frac{\sin (\pi\alpha\rho)}{\pi}
  x^{\alpha-1} , & x < 0,
  \end{IEEEcases}
\]
is invariant for the stable process killed upon hitting zero, and defines
the family of measures $(\stPaz_x)_{x \ne 0}$ given via the Doob $h$-transform:
\[
  \stPaz_x(\Lambda) = \frac{1}{h(x)} \stE_x[h(X_t) \Ind_\Lambda ; t < T_0],
  \wh x \ne 0,
\]
for $\Lambda \in \FF_t = \sigma(X_s, \, s \le t)$.
In \cite{Pan-cond} it is also shown that the laws $\stPaz_x$ arise as
limits of the stable process conditioned not to have hit zero up to
an exponential time of rate $q$, as $q \downto 0$.
The canonical process associated with the laws $(\stPaz_x)_{x \ne 0}$ is
therefore called the \define{stable process conditioned to avoid zero},
and we shall denote it by $\Xaz$.

\skippar
\newcommand{\LEaz}{\psi^{\updownarrow}}
Consider now the process $\Raz = \frac{1}{2} \abs{ \Xaz}$. This is
a pssMp, and we may consider its Lamperti transform, which we will
denote by $\xiaz$. The characteristics of the generalised Lamperti
representation of $\Xaz$ have been computed explicitly in \cite{CPR},
and the
Laplace exponent, $\LEaz$,
% characteristic exponent, $\CEaz$,
of $\xiaz$ could be
computed from this information; however, the harmonic transform
gives us the following straightforward relationship
between
Laplace exponents:
\[ \LEaz(z) = \LERp(z + \alpha - 1) . \]
% characteristic exponents:
% \[ \CEaz(\theta) = \CEp(\theta - \iu(\alpha-1)) .  \]
This allows us to calculate
\[
  \LEaz(2z)
  = -\frac{\Gamma(1/2-z)}{\Gamma((1-\alpha)/2-z)}
  \frac{\Gamma(\alpha/2+z)}{\Gamma(z)},
\]
% \[ \CEaz(2\theta)
%   = \frac{\Gamma(1/2-\iu\theta)}{\Gamma((1-\alpha)/2-\iu\theta)}
%   \frac{\Gamma(\alpha/2 + \iu\theta)}{\Gamma(\iu\theta)} ,
% \]
which demonstrates that $2\xiaz$ is a process in the extended
hypergeometric class, with parameters
\[ (\beta,\gamma,\betah,\gammah)
  = ((\alpha+1)/2, \alpha/2, 0, \alpha/2) . \]

The present authors and A.\ Kuznetsov previously computed $\LEaz$ in \cite{KKPW-T0},
where we also observed that the process $\xiaz$ is the dual L\'evy process to $\ltRp$,
and remarked that this implies a certain time-reversal relation between $R$ and $\Raz$;
see \cite[\S 2]{CP-lower}.

\section{Concluding remarks}

\newcommand\lsal{\boldsymbol{\alpha}}
\newcommand\lsbe{\boldsymbol{\beta}}
\newcommand\lsde{\boldsymbol{\delta}}
\newcommand\lscp{\mathbf{c}_{\boldsymbol{+}}}
\newcommand\lscm{\mathbf{c}_{\boldsymbol{-}}}
\newcommand\lsld{\boldsymbol{\pi}}

In this section, we offer some comments on how our approach may be adapted
in order to offer new insight on an existing class of processes, the Lamperti-stable
processes. These were defined in general in the work of \citet{CPP-LS},
and the one-dimensional Lamperti-stable processes are defined as follows.
We say that a L\'evy process $\xi$ is in the
\define{Lamperti-stable class} if it has no Gaussian component
and its L\'evy measure has density
\[
  \lsld(x)
  = \begin{cases}
      \lscp e^{\lsbe x} (e^x-1)^{-(\lsal+1)} \, \dd x, & x > 0, \\
      \lscm e^{-\lsde x} (e^{-x}-1)^{-(\lsal+1)} \, \dd x, & x < 0,
    \end{cases}
\]
for some choice of parameters
$(\lsal,\lsbe,\lsde,\lscp,\lscm)$ such that $\lsal\in(0,2)$
and $\lsbe,\lsde,\lscp,\lscm \ge 0$.

The one-dimensional Lamperti-stable processes
form a proper subclass of the $\beta$-class of L\'evy processes
of \citet{Kuz-beta}.
It was observed in \cite{KP-HG} that there is an intersection between
the hypergeometric class and the Lamperti-stable class. In particular,
the Lamperti representations of killed and conditioned stable processes
(see \cite{CC-LS}) fall within the hypergeometric class; and generally
speaking, setting $\beta=\betah$ in the hypergeometric class and choosing
$\gamma,\gammah$ as desired, one obtains a Lamperti-stable process.

However, not all Lamperti-stable processes may be obtained in this way,
and we now outline how the ideas developed in this work can be used to
characterise another subset of the Lamperti-stable processes.

\skippar
\newcommand{\AHL}{\mathcal{A}_{\text{\textnormal{EHL}}}}%
Define the set of parameters
\[ \AHL = \bigl\{ \beta \in [1,2], \, \gamma \in (1,2), \, \gammah \in (-1,0) \bigr\} \]
and for $(\beta,\gamma,\gammah) \in \AHL$, let
\[ \psi(z)
  = \frac{\Gamma(1-\beta+\gamma-z)}{\Gamma(2-\beta-z)}
  \frac{\Gamma(\beta+\gammah+z)}{\Gamma(\beta+z)}.
\]
Note that this is the negative of the usual hypergeometric Laplace exponent,
with $\beta=\hat\beta$. We claim that the following proposition holds.

\begin{prop}
  There exists a L\'evy process $\xi$ with Laplace exponent $\psi$.
  Its Wiener--Hopf factorisation
  $\psi(z)=-\kappa(-z)\kappah(z)$ is given by the components
  \[
    \kappa(z) = \frac{\Gamma(1-\beta+\gamma+z)}{\Gamma(2-\beta+z)},
    \qquad
    \kappah(z) = (\beta-1+z) \frac{\Gamma(\beta+\gammah+z)}{\Gamma(\beta+z)}.
  \]
  The ascending ladder height process is a Lamperti-stable subordinator,
  and the descending factor satisfies
  \[ \kappah(z) = \bigl(\Ttrans_{\beta-1}\upsilon\bigr)^*(z),
  \wh \upsilon(z) = \frac{\Gamma(1+z)}{\Gamma(1+\gammah+z)}. \]
  Here, $\upsilon$ is the Laplace exponent of a Lamperti-stable subordinator.

  The process $\xi$ has no Gaussian component and has a L\'evy density given by
  \[
    \pi(x) =
    \begin{cases}
      \dfrac{\Gamma(\gamma+\gammah+1)}{\Gamma(1+\gamma)\Gamma(-\gamma)}
      e^{(\beta+\gammah)x} (e^x-1)^{-(\gamma+\gammah+1)},
      & x > 0, \\
      \dfrac{\Gamma(\gamma+\gammah+1)}{\Gamma(1+\gammah)\Gamma(-\gammah)}
      e^{-(1-\beta+\gamma)x}
      (e^{-x}-1)^{-(\gamma+\gammah+1)},
      & x < 0. 
    \end{cases}
  \]
  Thus, $\xi$ falls within the Lamperti-stable class, and
  \[(\lsal,\lsbe,\lsde) = (\gamma+\gammah, \, \beta+\gammah, \, 1-\beta+\gamma) . \]
\end{prop}

The proposition may be proved in much the same way as \autoref{t:WHF XHG},
first using the theory of philanthropy to prove existence, and then
the theory of meromorphic L\'evy processes to deduce the L\'evy measure.

We have thus provided an explicit spatial Wiener--Hopf factorisation of a subclass of
Lamperti-stable processes disjoint from that given by the hypergeometric
processes. 

\paragraph{Acknowledgements} We would like to thank the anonymous referee for their thorough reading and
helpful comments, which have improved the article.

%%%%%%%%%%%%%%%%%%%%%%% generated by bibtex


\begin{thebibliography}{37}
\providecommand{\natexlab}[1]{#1}
\providecommand{\url}[1]{\texttt{#1}}
\expandafter\ifx\csname urlstyle\endcsname\relax
  \providecommand{\doi}[1]{doi: #1}\else
  \providecommand{\doi}{doi: \begingroup \urlstyle{rm}\Url}\fi

\bibitem[Bertoin(1996)]{Ber-Levy}
J.~Bertoin.
\newblock \emph{L\'evy processes}, volume 121 of \emph{Cambridge Tracts in
  Mathematics}.
\newblock Cambridge University Press, Cambridge, 1996.
\newblock ISBN 0-521-56243-0.

\bibitem[Bertoin and Yor(2002)]{BY-ent}
J.~Bertoin and M.~Yor.
\newblock The entrance laws of self-similar {M}arkov processes and exponential
  functionals of {L}\'evy processes.
\newblock \emph{Potential Anal.}, 17\penalty0 (4):\penalty0 389--400, 2002.
\newblock ISSN 0926-2601.
\newblock \doi{10.1023/A:1016377720516}.

\bibitem[Bertoin and Yor(2005)]{BY-ef}
J.~Bertoin and M.~Yor.
\newblock Exponential functionals of {L}\'evy processes.
\newblock \emph{Probab. Surv.}, 2:\penalty0 191--212, 2005.
\newblock ISSN 1549-5787.
\newblock \doi{10.1214/154957805100000122}.

\bibitem[Blumenthal and Getoor(1968)]{BG-mppt}
R.~M. Blumenthal and R.~K. Getoor.
\newblock \emph{Markov processes and potential theory}.
\newblock Pure and Applied Mathematics, Vol. 29. Academic Press, New York,
  1968.

\bibitem[Caballero and Chaumont(2006)]{CC-LS}
M.~E. Caballero and L.~Chaumont.
\newblock Conditioned stable {L}\'evy processes and the {L}amperti
  representation.
\newblock \emph{J. Appl. Probab.}, 43\penalty0 (4):\penalty0 967--983, 2006.
\newblock ISSN 0021-9002.
\newblock \doi{10.1239/jap/1165505201}.

\bibitem[Caballero et~al.(2010)Caballero, Pardo, and P{\'e}rez]{CPP-LS}
M.~E. Caballero, J.~C. Pardo, and J.~L. P{\'e}rez.
\newblock On {L}amperti stable processes.
\newblock \emph{Probab. Math. Statist.}, 30\penalty0 (1):\penalty0 1--28, 2010.
\newblock ISSN 0208-4147.

\bibitem[Caballero et~al.(2011)Caballero, Pardo, and P{\'e}rez]{CPP-radial}
M.~E. Caballero, J.~C. Pardo, and J.~L. P{\'e}rez.
\newblock Explicit identities for {L}\'evy processes associated to symmetric
  stable processes.
\newblock \emph{Bernoulli}, 17\penalty0 (1):\penalty0 34--59, 2011.
\newblock ISSN 1350-7265.
\newblock \doi{10.3150/10-BEJ275}.

\bibitem[Chaumont and Pardo(2006)]{CP-lower}
L.~Chaumont and J.~C. Pardo.
\newblock The lower envelope of positive self-similar {M}arkov processes.
\newblock \emph{Electron. J. Probab.}, 11:\penalty0 no. 49, 1321--1341, 2006.
\newblock ISSN 1083-6489.
\newblock \doi{10.1214/EJP.v11-382}.

\bibitem[Chaumont et~al.(2009)Chaumont, Kyprianou, and Pardo]{CKP-explicit}
L.~Chaumont, A.~E. Kyprianou, and J.~C. Pardo.
\newblock Some explicit identities associated with positive self-similar
  {M}arkov processes.
\newblock \emph{Stochastic Process. Appl.}, 119\penalty0 (3):\penalty0
  980--1000, 2009.
\newblock ISSN 0304-4149.
\newblock \doi{10.1016/j.spa.2008.05.001}.

\bibitem[Chaumont et~al.(2013)Chaumont, Pant{\'{\i}}, and Rivero]{CPR}
L.~Chaumont, H.~Pant{\'{\i}}, and V.~Rivero.
\newblock The {L}amperti representation of real-valued self-similar {M}arkov
  processes.
\newblock \emph{Bernoulli}, 19\penalty0 (5B):\penalty0 2494--2523, 2013.
\newblock ISSN 1350-7265.
\newblock \doi{10.3150/12-BEJ460}.

\bibitem[Cordero(2010)]{Cor-thesis}
F.~Cordero.
\newblock \emph{On the excursion theory for the symmetric stable L\'evy
  processes with index $\alpha \in ]1,2]$ and some applications}.
\newblock PhD thesis, Universit\'e Pierre et Marie Curie -- Paris VI, 2010.

\bibitem[Doney(1987)]{Don-Ckl}
R.~A. Doney.
\newblock On {W}iener-{H}opf factorisation and the distribution of extrema for
  certain stable processes.
\newblock \emph{Ann. Probab.}, 15\penalty0 (4):\penalty0 1352--1362, 1987.
\newblock ISSN 0091-1798.

\bibitem[Gnedin(2010)]{Gne-rrcs}
A.~V. Gnedin.
\newblock Regeneration in random combinatorial structures.
\newblock \emph{Probab. Surv.}, 7:\penalty0 105--156, 2010.
\newblock ISSN 1549-5787.
\newblock \doi{10.1214/10-PS163}.

\bibitem[Gradshteyn and Ryzhik(2007)]{GR}
I.~S. Gradshteyn and I.~M. Ryzhik.
\newblock \emph{Table of integrals, series, and products}.
\newblock Elsevier/Academic Press, Amsterdam, seventh edition, 2007.
\newblock ISBN 978-0-12-373637-6.
\newblock Translated from the Russian. Translation edited and with a preface by
  Alan Jeffrey and Daniel Zwillinger.

\bibitem[Kuznetsov(2010{\natexlab{a}})]{Kuz-beta}
A.~Kuznetsov.
\newblock Wiener-{H}opf factorization and distribution of extrema for a family
  of {L}\'evy processes.
\newblock \emph{Ann. Appl. Probab.}, 20\penalty0 (5):\penalty0 1801--1830,
  2010{\natexlab{a}}.
\newblock ISSN 1050-5164.
\newblock \doi{10.1214/09-AAP673}.

\bibitem[Kuznetsov(2010{\natexlab{b}})]{Kuz-theta}
A.~Kuznetsov.
\newblock Wiener-{H}opf factorization for a family of {L}\'evy processes
  related to theta functions.
\newblock \emph{J. Appl. Probab.}, 47\penalty0 (4):\penalty0 1023--1033,
  2010{\natexlab{b}}.
\newblock ISSN 0021-9002.

\bibitem[Kuznetsov and Pardo(2013)]{KP-HG}
A.~Kuznetsov and J.~C. Pardo.
\newblock Fluctuations of stable processes and exponential functionals of
  hypergeometric {L}\'evy processes.
\newblock \emph{Acta Appl. Math.}, 123:\penalty0 113--139, 2013.
\newblock ISSN 0167-8019.
\newblock \doi{10.1007/s10440-012-9718-y}.

\bibitem[Kuznetsov et~al.(2011)Kuznetsov, Kyprianou, Pardo, and van
  Schaik]{KKPvS}
A.~Kuznetsov, A.~E. Kyprianou, J.~C. Pardo, and K.~van Schaik.
\newblock A {W}iener-{H}opf {M}onte {C}arlo simulation technique for {L}\'evy
  processes.
\newblock \emph{Ann. Appl. Probab.}, 21\penalty0 (6):\penalty0 2171--2190,
  2011.
\newblock ISSN 1050-5164.
\newblock \doi{10.1214/10-AAP746}.

\bibitem[Kuznetsov et~al.(2012)Kuznetsov, Kyprianou, and Pardo]{KKP-mero}
A.~Kuznetsov, A.~E. Kyprianou, and J.~C. Pardo.
\newblock Meromorphic {L}\'evy processes and their fluctuation identities.
\newblock \emph{Ann. Appl. Probab.}, 22\penalty0 (3):\penalty0 1101--1135,
  2012.
\newblock \doi{10.1214/11-AAP787}.

\bibitem[Kuznetsov et~al.(2013)Kuznetsov, Kyprianou, and Rivero]{KKR-scale}
A.~Kuznetsov, A.~E. Kyprianou, and V.~Rivero.
\newblock The theory of scale functions for spectrally negative {L}\'evy
  processes.
\newblock In \emph{{L}\'evy {M}atters {II}}, Lecture Notes in Mathematics,
  pages 97--186. Springer, Berlin and Heidelberg, 2013.
\newblock ISBN 978-3-642-31406-3.
\newblock \doi{10.1007/978-3-642-31407-0_2}.

\bibitem[Kuznetsov et~al.(2014)Kuznetsov, Kyprianou, Pardo, and
  Watson]{KKPW-T0}
A.~Kuznetsov, A.~E. Kyprianou, J.~C. Pardo, and A.~R. Watson.
\newblock The hitting time of zero for a stable process.
\newblock \emph{Electron. J. Probab.}, 19:\penalty0 no. 30, 1--26, 2014.
\newblock ISSN 1083-6489.
\newblock \doi{10.1214/EJP.v19-2647}.

\bibitem[Kyprianou(2006)]{Kyp}
A.~E. Kyprianou.
\newblock \emph{Introductory lectures on fluctuations of {L}\'evy processes
  with applications}.
\newblock Universitext. Springer-Verlag, Berlin, 2006.
\newblock ISBN 978-3-540-31342-7.

\bibitem[Kyprianou and Patie(2011)]{KP-CT}
A.~E. Kyprianou and P.~Patie.
\newblock A {C}iesielski-{T}aylor type identity for positive self-similar
  {M}arkov processes.
\newblock \emph{Ann. Inst. Henri Poincar\'e Probab. Stat.}, 47\penalty0
  (3):\penalty0 917--928, 2011.
\newblock ISSN 0246-0203.
\newblock \doi{10.1214/10-AIHP398}.

\bibitem[Kyprianou and Rivero(2008)]{KR-scale}
A.~E. Kyprianou and V.~Rivero.
\newblock Special, conjugate and complete scale functions for spectrally
  negative {L}\'evy processes.
\newblock \emph{Electron. J. Probab.}, 13:\penalty0 1672--1701, 2008.
\newblock ISSN 1083-6489.
\newblock \doi{10.1214/EJP.v13-567}.

\bibitem[Kyprianou et~al.(2010)Kyprianou, Pardo, and Rivero]{KPR-n-tuple}
A.~E. Kyprianou, J.~C. Pardo, and V.~Rivero.
\newblock Exact and asymptotic {$n$}-tuple laws at first and last passage.
\newblock \emph{Ann. Appl. Probab.}, 20\penalty0 (2):\penalty0 522--564, 2010.
\newblock ISSN 1050-5164.
\newblock \doi{10.1214/09-AAP626}.

\bibitem[Kyprianou et~al.(2014)Kyprianou, Pardo, and Watson]{KPW-cens}
A.~E. Kyprianou, J.~C. Pardo, and A.~R. Watson.
\newblock Hitting distributions of {$\alpha$}-stable processes via path
  censoring and self-similarity.
\newblock \emph{Ann. Probab.}, 42\penalty0 (1):\penalty0 398--430, 2014.
\newblock ISSN 0091-1798.
\newblock \doi{10.1214/12-AOP790}.

\bibitem[Lamperti(1972)]{Lam-ssLT}
J.~Lamperti.
\newblock Semi-stable {M}arkov processes. {I}.
\newblock \emph{Z. Wahrscheinlichkeitstheorie und Verw. Gebiete}, 22:\penalty0
  205--225, 1972.

\bibitem[Lukacs and Sz{\'a}sz(1952)]{LS-acf}
E.~Lukacs and O.~Sz{\'a}sz.
\newblock On analytic characteristic functions.
\newblock \emph{Pacific J. Math.}, 2:\penalty0 615--625, 1952.
\newblock ISSN 0030-8730.

\bibitem[Olver(1974)]{Olv-asym}
F.~W.~J. Olver.
\newblock \emph{Asymptotics and special functions}.
\newblock Academic Press, New York--London, 1974.
\newblock Computer Science and Applied Mathematics.

\bibitem[Pant\'\i(2012)]{Pan-cond}
H.~Pant\'\i.
\newblock On {L}\'evy processes conditioned to avoid zero.
\newblock \arxivref{1304.3191v1}{math.PR}, 2012.

\bibitem[Pardo(2009)]{Par-upper}
J.~C. Pardo.
\newblock The upper envelope of positive self-similar {M}arkov processes.
\newblock \emph{J. Theoret. Probab.}, 22\penalty0 (2):\penalty0 514--542, 2009.
\newblock ISSN 0894-9840.
\newblock \doi{10.1007/s10959-008-0152-z}.

\bibitem[Rivero(2007)]{Riv-re2}
V.~Rivero.
\newblock Recurrent extensions of self-similar {M}arkov processes and
  {C}ram\'er's condition. {II}.
\newblock \emph{Bernoulli}, 13\penalty0 (4):\penalty0 1053--1070, 2007.
\newblock ISSN 1350-7265.
\newblock \doi{10.3150/07-BEJ6082}.

\bibitem[Rogozin(1972)]{Rog-hit}
B.~A. Rogozin.
\newblock The distribution of the first hit for stable and asymptotically
  stable walks on an interval.
\newblock \emph{Theory of Probability and its Applications}, 17\penalty0
  (2):\penalty0 332--338, 1972.
\newblock ISSN 0040585X.
\newblock \doi{10.1137/1117035}.

\bibitem[Song and Vondra{\v{c}}ek(2006)]{SV-special}
R.~Song and Z.~Vondra{\v{c}}ek.
\newblock Potential theory of special subordinators and subordinate killed
  stable processes.
\newblock \emph{J. Theoret. Probab.}, 19\penalty0 (4):\penalty0 817--847, 2006.
\newblock ISSN 0894-9840.
\newblock \doi{10.1007/s10959-006-0045-y}.

\bibitem[Vigon(2002)]{Vig-thesis}
V.~Vigon.
\newblock \emph{Simplifiez vos {L}\'evy en titillant la factorisation de
  {W}iener--{H}opf}.
\newblock PhD thesis, INSA de Rouen, 2002.

\bibitem[Vuolle-Apiala(1994)]{VA-Ito}
J.~Vuolle-Apiala.
\newblock It\^o excursion theory for self-similar {M}arkov processes.
\newblock \emph{Ann. Probab.}, 22\penalty0 (2):\penalty0 546--565, 1994.
\newblock ISSN 0091-1798.

\bibitem[Yano et~al.(2009)Yano, Yano, and Yor]{YYY-laws}
K.~Yano, Y.~Yano, and M.~Yor.
\newblock On the laws of first hitting times of points for one-dimensional
  symmetric stable {L}\'evy processes.
\newblock In \emph{S\'eminaire de {P}robabilit\'es {XLII}}, volume 1979 of
  \emph{Lecture Notes in Math.}, pages 187--227. Springer, Berlin, 2009.
\newblock \doi{10.1007/978-3-642-01763-6_8}.

\end{thebibliography}
\end{document}